\documentclass[11pt]{article}
\usepackage{amsfonts,amsmath,amscd,amssymb,amsthm}

\newfont{\eufm}{eufm10 scaled\magstep1}
\newcommand{\FRAK}[1]{\mbox{\eufm#1}}

\newcommand{\cO}{\mathcal{O}}

\newcommand{\cD}{\mathcal{D}}
\newcommand{\cB}{\mathcal{B}}
\newcommand{\cA}{\mathcal{A}}

\newcommand{\cR}{\mathcal{R}}

\newcommand{\bbN}{\mathbb{N}}
\newcommand{\bbZ}{\mathbb{Z}}
\newcommand{\bbC}{\mathbb{C}}

\newcommand{\bbQ}{\mathbb{Q}}

\newcommand{\bbF}{\mathbb{F}}

\newcommand{\bbX}{\mathbb{X}}
\newcommand{\bbY}{\mathbb{Y}}

\newtheorem{thm}{Theorem}[section]
\newtheorem{lem}[thm]{Lemma}

\newtheorem{prop}[thm]{Proposition}

\newtheorem{rem}[thm]{Remarks}

\begin{document}
\title{Finite dimensional representations of
invariant differential operators}
\author{Ian M. Musson\thanks{partially supported by NSF grant  DMS-0099923.}
\\ Department of Mathematical Sciences\\
University of Wisconsin-Milwaukee, USA\\E-mail: musson@uwm.edu\\\;\\
Sonia L. Rueda\\Departamento de \' Algebra y An\' alisis Matem\' atico \\
Universidad de Almer\'\i a\\04120 Almer\'\i a,
Spain\\E-mail:srueda@ual.es} \maketitle

\begin{abstract}
Let $k$ be an algebraically closed field of characteristic $0$,
$Y=k^{r}\times {(k^{\times})}^{s}$  and let $G$ be an algebraic
torus acting diagonally on the ring of differential operators $\cD
(Y)^G$. We give necessary and sufficient conditions for $\cD
(Y)^G$ to have enough simple finite dimensional representations,
in the sense that the intersection of the kernels of all the
simple finite dimensional representations is zero. As an
application we show that if $K\longrightarrow GL(V)$ is a
representation of a reductive group $K$ and if zero is not a
weight of a maximal torus of $K$ on $V$, then $\cD (V)^K$ has
enough finite dimensional representations. We also construct
examples of FCR- algebras with any GK dimension $\geq 3$.
\end{abstract}


\section{Introduction}\label{intro}
For a variety $Y$ over $k$ we denote the ring of regular functions
on $Y$ by $\cO (Y)$ and the ring of differential operators by $\cD
(Y)$. Recently there has been much interest in the study of the
invariant ring $\cD (Y)^G$ when $G$ is a reductive group acting on
a smooth affine variety $Y$, see for example \cite{LS}, \cite{M},
\cite{MV}, \cite{SG}, \cite{Sch}, \cite{VB}. In this paper our
primary focus is on the case where $Y=V\times W$ for a vector
space $V$, a torus $W$ and $G$ is a torus acting diagonally on
$Y$.

We say that a $k$-algebra $R$ has enough finite dimensional
modules ( resp. enough simple finite dimensional modules) if $\cap
ann_{R} M =0$ where the intersection is taken over all finite
dimensional (resp. simple finite dimensional) $R$-modules.

\vspace{0.5cm}{\bf Proposition A}
\begin{em}
If $\cD (Y)^G$ has a finite dimensional module then $G$ acts
transitively on $W$.
\end{em}
\vspace{0.5cm}

Now assume that $G$ acts transitively on $W$ and let $H$ be the
stabilizer in $G$ of $w\in W$. Note that the connected component
$H^o$ of the identity in $H$ is a torus, but we may have $H\neq
H^{o}$. The slice representation at $w$ is isomorphic to $(H,V)$,
see \S ~\ref{sr}. We give necessary and sufficient conditions for
$\cD(Y)^G$ to have enough finite dimensional simple modules.

\vspace{0.5cm} {\bf Theorem B}
\begin{em}
Assume that $G$ acts transitively on $W$. The following conditions
are equivalent.
\begin{enumerate}
\item $V^{H^o}=0$. \item $\cD (V)^{H^o}$ has enough simple finite
dimensional modules. \item $\cD (Y)^G$ has a nonzero finite
dimensional module. \item $\cD (Y)^G$ has enough simple finite
dimensional modules.
\end{enumerate}
\end{em}
\vspace{0.5cm}

Set $\FRAK{h}=Lie(H)\subseteq \FRAK{g}=Lie(G)$. For
$\lambda\in\FRAK{g}^*$, $\mu\in\FRAK{h}^*$ we set
\begin{equation}
\cB_{\lambda}(Y)=\cD (Y)^G/(\FRAK{g}-\lambda (\FRAK{g}))\;,\indent
\cB_{\mu}(V)=\cD (V)^H/(\FRAK{h}-\mu (\FRAK{h})).
\end{equation}
The algebra $B_{\lambda } (Y)$ is studied in detail in \cite{MV}.
Let $i^*:\FRAK{g}^*\longrightarrow\FRAK{h}^*$ be the map obtained
from the inclusion $i:\FRAK{h}\longrightarrow\FRAK{g}$.

\vspace{0.5cm}{\bf Proposition C}
\begin{em}
\begin{enumerate}
\item There is an injective algebra homomorphism $\xi:\cD
(V)^H\longrightarrow\cD (Y)^G$. \item If $\lambda\in\FRAK{g}^*$
and $\mu=i^*(\lambda )$ the above map induces an isomorphism
$\cB_{\mu}(V)\cong\cB_{\lambda} (Y)$.
\end{enumerate}
\end{em}
\vspace{0.5cm}

Note that any simple $\cD (Y)^G$-module is a $\cB_{\lambda}
(Y)$-module for some $\lambda\in\FRAK{g}^*$. So Propositions A and
C reduce the study of finite dimensional simple $\cD
(Y)^G$-modules to that of finite dimensional simple $\cD
(V)^H$-modules.  Some other situations in which slice
representations have been used to study invariant differential
operators may be found in \cite{SG} and \cite{VB}.

Now suppose that $G$ is a maximal torus in a reductive subgroup
$K$ of GL$(V)$. Then $\cD (V)^K$ is a subring of $\cD (V)^G$, so
any $\cD (V)^G$-module is also a $\cD (V)^K$-module. Thus we
obtain as an immediate corollary to Theorem B.

\vspace{0.5cm}{\bf Theorem D}
\begin{em}
If $0$ is not a weight of $G$ on $V$ then $\cD (V)^K$ has enough
finite dimensional modules.
\end{em}
\vspace{0.5cm}

In general it seems rather difficult to say much about finite
dimensional $\cD (V)^K$-modules (or more generally primitive
ideals in $\cD (V)^K$) if $K^{o}$ is not a torus. As far as we are
aware if $K$ is simple the only cases where anything is known are
the case of the adjoint representation of $SL(3)$ \cite{Sch} and
some representations arising from classical invariant theory
\cite{LS}.

Here is a brief outline of the paper. First we establish some
notation and express the conditions in Proposition A and Theorem B
in terms of the weights of the torus action. If either the action
of $G$ on $W$ is not transitive, or if $V^{H^o}\neq 0$ we show
that $\cD (Y)^G$ contains the fixed ring of the first Weyl algebra
under the action of a finite group. If neither of these conditions
hold we show that $\cD (Y)^G$ has enough finite dimensional
modules by first analyzing the case $\cO (Y)^G=k$, and then using
Fourier transforms to reduce to this case.

Kraft and Small \cite{KS} call an algebra $R$ an {\em FCR-algebra}
if $R$ has enough finite dimensional modules and every finite
dimensional $R$-module is completely reducible. In the last
section of the paper we combine Theorem B with results from
\cite{MV} to give examples of FCR algebras with any given integer
Gelfand-Kirillov dimension $\geq 3$. We denote the
Gelfand-Kirillov dimension of an algebra $A$ by GK-dim$A$, see
\cite{KL} for background.

\section{Notation}\label{notation}
\subsection{Two exact sequences}
We assume that $V = k^{r}, W = (k^{\times})^{s}$, and $Y = V
\times W\subseteq k^n$, where $n=r+s$. A {\em diagonal action} of
a torus $G$ on $Y$ is an action that extends to a diagonal action
on $k^n$. Such an action is given by an embedding of $G$ into the
group $T$ of diagonal matrices in $GL(n)$. Set $\bbX
(G)=Hom(G,k^{\times})$, $\bbY
(G)=Hom(k^{\times},G)=Hom(\bbX(G),\bbZ)$, the groups of characters
and one-parameter subgroups of $G$, respectively. Let
\begin{align}
&M=\bbY (G),\; N=\bbX (G)\\
&A=\bbY (T),\; B=\bbX (T).
\end{align}
There are natural bilinear pairings
\begin{eqnarray}
<\;,\;>:M\times N\longrightarrow \bbZ\\
(\;\;,\;):A\times B\longrightarrow \bbZ.
\end{eqnarray}
Let $\phi :\bbX (T)\longrightarrow \bbX (G)$ be the restriction
map. Let $K=ker\phi$ and
\begin{equation}
K^{\perp}=\{\chi\in\bbY (T) | (\chi ,K)=0\}.
\end{equation}
Then we have exact sequences
\begin{align}
&0\longrightarrow K^{\perp}\longrightarrow \bbY
(T)=A\longrightarrow \bbY (T)/K^{\perp}\longrightarrow 0\\
&0\longrightarrow K\longrightarrow \bbX (T)=B\longrightarrow
N'\longrightarrow 0.
\end{align}
where $N'=\mbox{Image}\phi\subseteq N=\bbX (G)$. Finally there is
a homomorphism
\begin{equation}
\omega:M\longrightarrow K^{\perp}
\end{equation}
defined by
\begin{equation}\label{w}
(\omega(m),b)=<m,\phi (b)>
\end{equation}
for $m\in M$, $b\in B$. A slightly different situation is
described in \cite{M}.

By introducing bases we can describe these maps using matrices.
Identify $G$ with $(k^{\times})^m$ and for $1\leq i\leq m$ define
$\nu_i\in\bbX (G)$ by $\nu_i (g)=g_i$ for $g=(g_1,\ldots ,g_m )\in
G$. We call $\{\nu_i\}$ the {\em standard basis} for $N$. Identify
$T$ with $(k^{\times})^n$ the standard basis $\{e_i\}$ for $B$ is
defined similarly, and write $\{\nu_i^*\}$, $\{e_i^*\}$ for the
dual bases of $M$ and $A$. For matrix computations we identify
$A$, $B$, $M$, $N$ as $\bbZ$-modules of column vectors using these
bases. Let $I_t$ denote the $t\times t$ identity matrix. If $E$ is
an abelian group we set $E_{\bbQ}=E\otimes_{\bbZ}\bbQ$.

For $1\leq i\leq n$ let $\eta_i=\phi (e_i )$. There is an integer
matrix $L=(l_{ij})$ such that
\begin{equation}
\eta_i=\sum_{j=1}^m l_{ji} \nu_j
\end{equation}
for $i=1,\ldots ,n$. From ~\eqref{w} it follows that
\begin{equation}
\omega(\nu_i^*)=\sum_{j=1}^n l_{ij}e_j^*.
\end{equation}
Thus the maps $B\longrightarrow N$ and $M\longrightarrow A$ are
given by multiplication by $L$ and the transpose of $L$
respectively. If $\beta=(\beta_1,\ldots ,\beta_n)$ belongs to the
submodule $R$ of $\bbZ^n$ spanned by the rows of $L$ then
\begin{equation}
\beta  (t)=(t^{\beta_1}, t^{\beta_2},\ldots ,t^{\beta_n})\in T
\end{equation}
determines the action of a one-parameter subgroup of $G$ on $Y$.
We summarize this data by saying that {\em G acts on $Y$ by the
matrix $L$}.

\begin{lem}\label{matrix}
If $\eta_{r+1},\ldots ,\eta_n$ are linearly independent there
exist matrices $\Gamma\in GL_m (\bbZ)$, $\Delta\in GL_n (\bbZ)$
such that
\begin{enumerate}
\item $\Gamma L\Delta$ has the block matrix form
\begin{equation}
\Gamma
L\Delta=\left[\begin{array}{cc}L_1&0\\L_2&D\end{array}\right],
\end{equation}
where $D$ is a diagonal matrix with nonzero diagonal entries
$d_1,\ldots ,d_s$.

\item $\Delta$ has the block matrix form
\begin{equation}
\Delta=\left[\begin{array}{cc}I_r&0\\0&\Delta_1\end{array}\right]
\end{equation}
with $\Delta_1\in GL_s (\bbZ)$.
\end{enumerate}
\end{lem}
\begin{proof}
Let $R$ be the submodule of $\bbZ^n$ spanned by the rows of $L$,
and let $\epsilon :\bbZ^n\longrightarrow \bbZ^s$ be the projection
onto the last $s$ coordinates. Since the submatrix of $L$ obtained
by deleting the first $r$ columns has rank $s$ it follows that
$\epsilon (R)$ is free abelian of rank $s$. Thus if
$R'=ker\epsilon$, the map $R\longrightarrow R/R'$ splits and there
is a submodule $R''$ of $R$ such that $R=R'\oplus R''$. By
choosing bases for $R'$ and $R''$ we find $\Gamma_1\in GL_m (\bbZ
)$ such that
\begin{equation*}
\Gamma_1L=\left[\begin{array}{cc}L_1&0\\L'_2&L_3\end{array}\right].
\end{equation*}
where $L_1$ has $m-s$ rows and $r$ columns and the rows of $[L_1
,0]$ form a basis for $R'$. There exist matrices $\Gamma_2\in GL_s
(\bbZ )$ and $\Delta_1\in GL_s (\bbZ)$ such that
$D=\Gamma_2L_3\Delta_1$ has the desired form, and we set
\begin{equation*}
\Gamma=\left[\begin{array}{cc}I_{m-s}&0\\0&\Gamma_2\end{array}\right]\cdot
\Gamma_1,\indent
\Delta=\left[\begin{array}{cc}I_{r}&0\\0&\Gamma_1\end{array}\right].
\end{equation*}
\end{proof}

Now suppose $G$ acts on $Y$ via the matrix $L$ and set
\begin{equation}
\Sigma_L=\{\alpha\in\bbN^r\times\bbZ^s |L\alpha =0\}.
\end{equation}
If $\alpha=(\alpha_1,\ldots ,\alpha_n)\in\bbN^r\times\bbZ^s$, we
set $Q^{\alpha }=Q_{1}^{\alpha_1}\ldots Q_{n}^{\alpha_n}$. Then
$\cO (Y)^G=span\{Q^{\alpha}\in\cO (Y) |L\alpha=0\}\cong k\Sigma_L$
the semigroup algebra on $\Sigma_L$. If $\Gamma$, $\Delta$ are as
in the lemma and $L'=\Gamma L\Delta$, there is an isomorphism
\begin{equation}
\Sigma_L\longrightarrow \Sigma_{L'}
\end{equation}
given by $x\mapsto \Delta^{-1}x$. Thus if $\eta_{r+1},\ldots
,\eta_n$ are linearly independent we assume henceforth that $L$
has the special form
\begin{equation}\label{newL}
L=\left[\begin{array}{cc}L_1&0\\L_2&D\end{array}\right]
\end{equation}
where $D$ is a diagonal matrix with nonzero diagonal entries
$d_1,\ldots ,d_s$.

\subsection{Rings of differential operators}
Note that $Y$ is a toric variety with a dense torus
$T=(k^{\times})^n\subseteq Y$. Write $Q_i$ for the character $e_i$
considered as a regular function on $T$, and let $P_i=\partial /
\partial Q_i$. Then
\begin{equation}\label{OX}
\cO (Y)=k[Q_1,\ldots ,Q_{r},Q_{r+1}^{\pm 1},\ldots ,Q_n^{\pm
1}],\end{equation}
\begin{equation}\label{DX}
\cD(Y)=k[Q_1,\ldots ,Q_{r},Q_{r+1}^{\pm 1},\ldots ,Q_n^{\pm
1},P_1,\ldots ,P_n].
\end{equation}

We consider the action of $G$ on $O (T)$ (or $O (Y)$) given by
{\em right} translation, that is
\begin{equation}
(g\centerdot f)(t)=f(tg)
\end{equation}
for $g\in G$, $f\in \cO (T)$, $t\in T$. This convention implies
that $Q_i$ has weight $\eta_i$ as in \cite{MV}, section 6.

This action extends to an action of $G$ on $\cA =k[Q_1^{\pm
1},\ldots,Q_n^{\pm 1},P_1, \ldots ,P_n]$. It is often convenient
to use exponential notation for elements of $\cA$. If $\lambda
=(\lambda_1,\ldots ,\lambda_n )\in\bbZ^n$ and $\mu =(\mu_1,\ldots
,\mu_n)\in\bbN^n$, we set $Q^{\lambda }=Q_{1}^{\lambda_1}\ldots
Q_{n}^{\lambda_n}$, and $P^{\mu}=P_{1}^{\mu_{1}}\ldots
P_{n}^{\mu_{n}}$. We have for $g=(g_1,\ldots ,g_m )\in G$
\begin{equation}\label{act}
g\centerdot
Q^{\lambda}=\prod_{i=1}^{m}g_{i}^{\sum_{j=1}^{n}\lambda_{j}l_{ij}}
Q^{\lambda}, \indent g\centerdot
P^{\mu}=\prod_{i=1}^{m}g_{i}^{\sum_{j=1}^{n} -\mu_{j}l_{ij}}
P^{\mu}.
\end{equation}
The elements $Q^{\lambda} P^{\mu} \in \cD (Y)$ with
$L\lambda=L\mu$ form a basis of $\cD (Y)^G$.

For $\alpha\in\bbZ^n$ we define
\begin{equation}\label{u}
u_{\alpha}=Q_1^{(\alpha_1)}\cdots Q_{r}^{(\alpha_r)}
Q_{r+1}^{\alpha_{r+1}}\cdots Q_n^{\alpha_n}
\end{equation}
where
\begin{equation}\label{x1}
Q_i^{(\alpha_i)}=\left\{\begin{array}{c}
Q_i^{\alpha_i}\mbox{ if }\alpha_i\geq 0\\
P_i^{-\alpha_i}\mbox{ if }\alpha_i< 0
\end{array}\right.,
\indent i\in\{1,\ldots ,r\}.
\end{equation}

If $\Pi_i=Q_iP_i$, then $[\Pi_i, u_{\alpha}]=\alpha_i u_{\alpha}$.
Therefore $\cD (Y)=\oplus\cD (Y)_{\alpha}$ where $\cD
(Y)_{\alpha}=\cA_0 u_{\alpha}$, $\alpha\in\bbZ^n$, with
$\cA_0=k[\Pi_1,\ldots ,\Pi_n]$. Then
\begin{equation}\label{AG}
\cD(Y)^G=\oplus_{\alpha\in Supp\cD(Y)^G} \cD(Y)_{\alpha}
\end{equation}
with $Supp \cD (Y)^G=\{\alpha\in\bbZ^n | L\alpha =0\}$.

Define $\Lambda\subseteq \bbZ^m$ by $\Lambda=\{L\alpha
|\alpha\in\bbN^r\times\bbZ^s\}$. For $\chi \in \Lambda$ define
\begin{equation}
\cO (Y)_{\chi} =span\{ Q^{\lambda}\in \cO (Y) | L\lambda=\chi \}.
\end{equation}
It is easy to see that
\begin{equation}
\cO (Y)=\oplus_{\chi\in \Lambda} \cO (Y)_{\chi}.
\end{equation}

Each $\cO (Y)_{\chi}$ is a $\cD (Y)^{G}$-module: for a given
$Q^{\lambda}\in\cO (Y)_{\chi}$ and $Q^{\mu}P^{\nu}\in \cD (Y)^G$,
$L\lambda =\chi$ and  $L\mu =L\nu$. Hence $Q^{\mu}
P^{\nu}\centerdot Q^{\lambda}$ is a multiple of $Q^{\lambda + \mu
-\nu}$ with $L(\lambda+ \mu -\nu)= L\lambda=\chi$.

\subsection{Some remarks on the algebra $\cD (V)^H$}\label{rem}

Note that $\cD (V)^{H}$ is the fixed ring of $\cD (V)^{H^o}$ under
the action of the finite abelian group $H/H^o$. Thus $\cD
(V)^{H^o}$ is graded by the character group of $H/H^o$ and $\cD
(V)^H$ is the identity component of this grading. Hence if $I$ is
any right ideal of $\cD (V)^H$ we have
\begin{equation}
I\cD (V)^{H^o}\cap\cD (V)^H=I.
\end{equation}
Now fix $\mu\in\FRAK{h}^*$ and let \FRAK{m} be the ideal of
$\cA'_0=k[\Pi_1,\ldots ,\Pi_r]$ generated by
$\FRAK{h}-\mu(\FRAK{h})$. We can apply the above remarks with
$I=\FRAK{m}\cD (V)^H$. It follows that there is an embedding of
algebras
\begin{equation}
\frac{\cD (V)^H}{\FRAK{m}\cD (V)^H}\longrightarrow\frac{\cD
(V)^{H^o}}{\FRAK{m}\cD (V)^{H^o}}.
\end{equation}
From now on we denote these algebras simply by $\cD
(V)^H/(\FRAK{h}-\mu (\FRAK{h}))$ and $\cD (V)^{H^o}/(\FRAK{h}-\mu
(\FRAK{h}))$. Note that by \cite{KL}, Lemma 3.1 and \cite{MV},
Theorem 8.2.1 we have
\begin{equation}
GK-dim \frac{\cD (V)^H}{(\FRAK{h}-\mu (\FRAK{h}))}\leq 2(n-m).
\end{equation}
In fact it is easily seen that equality holds.

\section{Actions of tori and slice representations}\label{sr}

In this section we express the hypothesis in Proposition A and
Theorem B in terms of the weights of the action.

\begin{lem}\label{li}
Suppose the torus $G=(k^\times )^m$ acts on $W=(k^{\times})^s$
with weights $\eta_{r+1},\ldots ,\eta_n$ and $w=(w_{r+1},\ldots
,w_n)\in W$.
\begin{enumerate}
\item If $\eta_{r+1},\ldots ,\eta_n$ are linearly independent then
$G$ acts transitively on $W$. \item If $\eta_{r+1},\ldots ,\eta_n$
are linearly dependent then any orbit of $G$ on $W$ has dimension
less than dim$W$.
\end{enumerate}
\end{lem}
\begin{proof}
\begin{enumerate}
\item Suppose that $\eta_{r+1},\ldots ,\eta_n $ are linearly
independent and let $L'=(l'_{ji})$ be the submatrix of $L$ spanned
by the last $s$ columns $\eta_{r+1},\ldots ,\eta_n$. Thus
$l'_{ji}=l_{j,i+r}$ for $1\leq j\leq m$ and $1\leq i\leq s$. Since
$L'$ has rank $s$ there is an integer matrix $F$ and a nonzero
integer $d$ such that $FL'=dI_s$.

To show that $G$ acts transitively on $W$ it suffices to show that
given $w=(w_{r+1},\ldots ,w_n)\in W$ there exists $g\in G$ such
that $g\cdot (1,\ldots ,1)=w$. For each $i$ choose $v_i\in
k^{\times}$ such that $v_i^d=w_i$ and set
$g_j=\prod_{p=1}^sv_p^{f_{p,j}}$ and $g=(g_1,\ldots ,g_m)\in G$.
Then $\eta_{i+r}(g)=\prod_{j=1}^m g_j^{l_{j,i+r}}=v_i^d=w_i$ as
required.

\item We can assume that $G$ acts faithfully on $W$. Then
$\cap_{i=r+1}^n ker\eta_i=1$. If $\eta_{r+1},\ldots ,\eta_n$ are
linearly dependent this implies that dim$G=m<s$. Any element $w\in
W$ has trivial stabilizer, so dim$Gw=$dim$G<$dim$W$.
\end{enumerate}
\end{proof}

Suppose that
\begin{equation}
\phi : G\times X\longrightarrow X
\end{equation}
is an action of a reductive group $G$ on a variety $X$ and
consider a point $u\in X$ with stabilizer $H=G_u$ and tangent
space $T_u (X)$. The differential
\begin{equation}
d\phi :\FRAK{g}\longrightarrow T_u (X)
\end{equation}
has kernel $\FRAK{h}=Lie(H)$. If $H$ is reductive we have
\begin{equation}
T_u (X)=d\phi (\FRAK{g}/\FRAK{h})\oplus U
\end{equation}
for some $H$-module $U$. We call the pair $(H,U)$ the slice
representation at $u$, see \cite{Lu}, \cite{Sl}.

Now suppose $G=(k^\times )^m$ a torus acting faithfully on
$Y=V\times W$ with weights $\eta_1,\ldots ,\eta_n$. Suppose that
$G$ acts transitively on $W$, and let $w=(w_{r+1},\ldots ,w_n)$ be
an element of $W$. Then
\begin{equation}\label{H}
H=G_w=\cap_{i=r+1}^n ker\eta_i.
\end{equation}

\begin{lem}
The slice representation at $w$ is isomorphic to $(H,V)$.
\end{lem}
\begin{proof}
Assume that $G$ acts on $Y$ via the matrix $L$ as in
~\eqref{newL}. Then we can identify $\FRAK{g}$ with
\begin{equation}
span\{\sum_{j=1}^n b_jQ_jP_j | b\in R\},
\end{equation}
with $R$ as in the proof of Lemma ~\ref{matrix}. The maximal ideal
$\FRAK{m}_w$ of $w$ in $\cO (Y)$ is generated by the $Q_i-w_i$
($1\leq i\leq n$) where $w_1=\ldots =w_{r}=0$. For $a\in \cO (Y)$
write $\overline{a}$ for the image of $a$ mod $\FRAK{m}_w^2$. For
$1\leq i\leq n$ let $T_i=\overline{Q_i-w_i}$ and let $\{T_i^*\}$
be the basis of $T_w (Y)$ dual to the basis $\{T_i\}_{1\leq i\leq
n}$ of $\FRAK{m}_w/\FRAK{m}_w^2$. The map $d\phi
:\FRAK{g}\longrightarrow T_w (Y)$ is defined by $x\cdot
\overline{a}= e_w (x\cdot a)$ for $x\in \FRAK{g}$ and $a\in
\FRAK{m}_w$ where $e_w: \cO (Y)\longrightarrow k$ is the
evaluation at $w$. An easy calculation shows that $Q_jP_j\cdot
T_i=\delta_{ij}w_i$, $d\phi$ maps $\FRAK{g}$ to the subspace of
$T_w (Y)$ spanned by $T_{r+1}^*,\ldots ,T_n^*$. The linear
independence of $\eta_{r+1},\ldots ,\eta_n$ implies that
dim$\FRAK{g}/\FRAK{h}=s$, so comparing dimensions we see that
$d\phi$ induces an isomorphism of $\FRAK{g}/\FRAK{h}$ onto
$span\{T_{r+1}^*,\ldots ,T_n^*\}$.
\end{proof}

Finally $H^o$ is the subtorus of $G$ generated by images of one
parameter subgroups corresponding to rows of $L_1$. For $1\leq
i\leq r$ let $\rho_i$ be the restriction of $\eta_i$ to $H$. These
characters can be thought of as columns of $L_1$.

We can easily see the following,
\begin{lem}\label{nz}
$V^{H^o}=0$ if and only if $\rho_i\neq 0$ for all $i=1,\ldots ,r$.
\end{lem}

\section{The special case $\cO (Y)^G=k$}

\begin{lem}\label{special}
\begin{enumerate}
\item If $\cO (Y)^G=k$ then $\eta_{r+1},\ldots ,\eta_n$ are
linearly independent.

\item Assume $\eta_{r+1},\ldots ,\eta_n$ are linearly independent.
Then the following conditions are equivalent.
\begin{enumerate}
\item $\cO (Y)^G=k$. \item For all $\chi\in\Lambda $, $\cO
(Y)_{\chi}$ is finite dimensional. \item There exists
$\beta=(\beta_1,\ldots ,\beta_n)\in K^{\perp}$ such that $\beta_i
>0$ for $i=1,\ldots ,r$, and $\beta_i=0$ for $i=r+1,\ldots ,n$.
\end{enumerate}
\end{enumerate}
\end{lem}
\begin{proof}
1. If $\eta_{r+1},\ldots ,\eta_n$ are not linearly independent we
can write
\begin{equation}
\sum_{i\in I}a_i\eta_i =\sum_{j\in J}b_j\eta_j
\end{equation}
where $I,J$ are disjoint subsets of $\{r+1,\ldots ,n\}$,
$I\neq\emptyset$ and the $a_i,b_j$ are positive integers. Then
\begin{equation}
\Pi_{j\in J} Q_j^{b_j}/ \Pi_{i\in I} Q_i^{a_i}
\end{equation}
is a nonconstant element of $\cO (Y)^G$.

2. (a)$\implies$(c) Suppose $\cO (Y)^G=k$. For $1\leq i\leq r$ set
\begin{equation}
C_i=\sum_{j=1,j\neq i}^r \bbQ_{\geq 0}\eta_j +\sum_{j=1}^s
\bbQ\eta_{r+j}\subseteq \bbX (G)_{\bbQ}
\end{equation}
and
\begin{equation}
C_i^{\vee}=\{\gamma\in\bbY (G)_{\bbQ} | <\gamma ,C_i>\geq 0\}.
\end{equation}
If $\alpha\in K$ with $\alpha_i\geq 0$ for $i=1,\ldots r$ we have
$\alpha=0$ since $Q^{\alpha}\in \cO (Y)^G$. It follows that
$-\eta_i\notin C_i$ for $1\leq i\leq r$. Hence by equation (*) on
page 9 of \cite{F}, there exists $\gamma_i\in C_i^{\vee}$ such
that $<\gamma_i ,\eta_i>=(\omega (\gamma_i),e_i)>0$. If
$\gamma=\sum_{i=1}^r\gamma_i$ then some integer multiple $\beta$
of $\omega (\gamma)$ satisfies the condition in (c).

(c)$\implies$(b) We may assume that $L$ has the form
~\eqref{newL}. Suppose that $\chi\in\Lambda$ and fix
$\varphi\in\bbN^r\times\bbZ^s$ with $L\varphi=\chi$. If $B$ is
identified with $\bbZ^n$ we have
\begin{equation}
\cO (Y)_{\chi}=span\{Q^{\alpha}|\alpha\in \bbN^r\times\bbZ^s ,
L(\alpha -\varphi )=0\}.
\end{equation}
For $\alpha\in B$ write $\alpha=\left(\begin{array}{c}\alpha' \\
\alpha'' \end{array}\right)$ with $\alpha' \in\bbN^r$ and
$\alpha''\in\bbZ^s$. Define $\varphi'$ and $\varphi''$ similarly.
To show $\cO (Y)_{\chi}$ is finite dimensional it suffices to show
there are at most finitely many $\alpha\in\bbN^r\times\bbZ^s$ such
that
\begin{equation}\label{L1}
L_1 (\alpha' -\varphi')=0
\end{equation}
\begin{equation}\label{L2}
 L_2(\alpha'-\varphi')+D(\alpha'' -\varphi'')=0.
\end{equation}
Since $D$ is invertible $\alpha''$ is determined once we fix
$\alpha'$. (Note however that if $D^{-1}$ does not have integer
entries, equation ~\eqref{L2} may impose additional conditions on
$\alpha'$).

Thus we may assume that $r=n$ and $L_1=L$. The condition $L
(\alpha -\varphi)=0$ is equivalent to $(K^{\perp},\alpha-\varphi
)=0$. Hence given $\beta\in K^{\perp}$ as in (c) we have
\begin{equation}
\sum_{i=1}^n\beta_i\alpha_i =\sum_{i=1}^n\beta_i\varphi_i
\end{equation}
and this equation has only finitely many solutions for
$\alpha\in\bbN^n$ since all $\beta_i$ are positive.

This completes the proof since obviously (b)$\implies$ (a).
\end{proof}

\begin{rem}
\begin{enumerate}
\item The first part of the Lemma can be used to give another
proof of one implication of Lemma ~\ref{li}.

\item If $k=\bbC$ and $Y=V$ there is a more geometric proof of the
equivalence of (a) and (c) in the Lemma. Indeed (a) says that $0$
is the only closed orbit of $G$ on $V$, while if $\beta$ is a
one-parameter subgroup of $G$ as in (c) then
$\displaystyle{\lim_{t\rightarrow 0}} \beta(t) v=0$ for all $v\in
V$. Thus the equivalence follows by the Hilbert-Mumford criterion
\cite{K}, III 2.2.
\end{enumerate}
\end{rem}

Let $Q_i^{(\alpha_i)}$ as in ~\eqref{x1}. The following identities
are easily proved:
\begin{equation}\label{mult1}
\mbox{ if }\alpha_i\geq 0\mbox{, then } Q_i^{(\alpha_i)}\centerdot
Q_i^{\lambda_i}=Q_i^{\lambda_i+\alpha_i}
\end{equation}
\begin{equation}\label{mult2}
\mbox{if }\alpha_i< 0\mbox{, then } Q_i^{(\alpha_i)}\centerdot
Q_i^{\lambda_i}=\left\{\begin{array}{l} \frac{\lambda_i
!}{(\lambda_i +\alpha_i)!}Q_i^{\lambda_i+\alpha_i}
\mbox{ if }-\alpha_i\leq\lambda_i\\
0\mbox{ if }-\alpha_i >\lambda_i
\end{array}\right. .
\end{equation}

\begin{lem}\label{simple}
If $\cO(Y)_{\chi }\neq 0$, then it is a simple $\cD
(Y)^{G}$-module.
\end{lem}
\begin{proof}
Let $Q^{\lambda},\;Q^{\mu}\in\cO (Y)_{\chi}$. Then $L\lambda
=L\mu=\chi$. Hence $u_{\mu -\lambda}\in \cD (Y)^{G}$. Set
$\alpha=\mu -\lambda$. If $1\leq i\leq r$, then $\mu_i\geq 0$, so
from ~\eqref{mult1}, ~\eqref{mult2} we get $u_{\mu -\lambda}
\centerdot Q^{\lambda}=c\; Q^{\mu }$, where $c$ is a nonzero
integer. Since all weight spaces of $\cO (Y)_{\chi}$ are one
dimensional the result follows from this.
\end{proof}

Assume that $\cO (Y)^G=k$. We assume the action of $G$ on $Y$ is
defined via the matrix $L$ in ~\eqref{newL}. The dimensions of the
modules $\cO (Y)_{\chi}$ can be calculated using the following
result. To state it we require some notation. Let $\bbF$ be a
product of cyclic groups of order $d_1,\ldots, d_s$ and
$R=\bbZ[\bbF][[t_1,\ldots ,t_{m-s}]]$ a ring of formal power
series over the group $\bbZ[\bbF]$. For $1\leq i\leq s$ let
$t_{m-s+i}$ be a generator for the cyclic subgroup of $\bbF$ of
order $d_i$. For $\chi=(\chi_1,\ldots ,\chi_m)\in\bbZ^m$ set
$t^{\chi}=t_1^{\chi_1} t_2^{\chi_2} \cdots t_m^{\chi_m}$.

\begin{prop}
The dimensions of the modules $\cO (Y)_{\chi}$ satisfy
\begin{equation}
\sum_{\chi\in\Lambda}dim\cO (Y)_{\chi} \; t^{\chi}=\prod_{j=1}^r
(1-\prod_{i=1}^m t_i^{l_{ij}})^{-1}.
\end{equation}
\end{prop}
\begin{proof}
The coefficient of $t^{\chi}$ in the expansion of the right side
equals the number of solutions for $\alpha\in\bbN^r$ to the
equations
\begin{equation}
\chi_i =\sum_{j=1}^r l_{ij}\alpha_j
\end{equation}
for $1\leq i\leq m-s$, and the congruences
\begin{equation}
\chi_{i+m-s}\equiv \sum_{j=1}^r l_{i+m-s,\; j}\alpha_j\;\mbox{
mod}\; d_i.
\end{equation}
for $1\leq i\leq s$.  This is easily seen to equal the number of
solutions to equations ~\eqref{L1} and ~\eqref{L2}.
\end{proof}

\section{Reduction to the special case}\label{csw}
Suppose $I\subseteq\{1,\ldots ,r\}$. For $1\leq i\leq n$ set
\begin{equation}\label{LI}
\varsigma_i=\left\{ \begin{array}{c}
-\eta_i\mbox{ if }i\in I\\
\eta_i\mbox{ if }i\notin I
\end{array}\right. .
\end{equation}
Let $L_I$ be the matrix with columns $\varsigma_1,\ldots
,\varsigma_n$. Then $L_I$ defines an action of a new torus $G_I$
on $Y$. As before $G_I$ can be thought of as a subtorus of $T$.

\begin{lem}\label{I}
If $V^{H^o}=0$, there is a subset $I$ of $\{1,\ldots ,r\}$ such
that $\cO (Y)^{G_I}=k$.
\end{lem}
\begin{proof}
Let $R'$ be the submodule of $\bbZ^{m-s}$ spanned by the rows of
$L_1$, and let $\epsilon_j :R_{\bbQ}'\longrightarrow \bbQ$ be the
restriction of the projection onto the $j^{th}$ coordinate. Since
$V^{H^o}= 0$ it follows that $\epsilon_j (R_{\bbQ}')\neq 0$ for
all $j$ so $ker\epsilon_j\neq R'_{\bbQ}$. Therefore $\cup_{j=1}^r
ker\epsilon_j\neq R'_{\bbQ}$ and we can find $\beta\in R'_{\bbQ}$
such that $\epsilon_j (\beta)\neq 0$ for $j=1,\ldots ,r$. Now set
$I=\{j\in\{1,\ldots ,r\}|\epsilon_j (\beta)<0\}$. Using Lemma
~\ref{special} it is easy to see that $\cO (Y)^{G_I}=k$.
\end{proof}

\begin{lem}\label{isomorphism}
Let $I$ be a subset of $\{1,\ldots ,r\}$. Then the map
$\sigma_I:\cA\rightarrow \cA$ defined by
\begin{equation}\label{sigmaI}
\sigma_I (Q_i)=\left\{\begin{array}{c}
-P_i\mbox{ if }i\in I\\
Q_i\mbox{ if }i\notin I
\end{array}\right.
\indent
\sigma_I (P_i)=\left\{ \begin{array}{c}
Q_i\mbox{ if }i\in I\\
P_i\mbox{ if }i\notin I
\end{array}\right.
\end{equation}
$i=1,\ldots ,n$ is an isomorphism between $\cD (Y)^G$ and $\cD
(Y)^{G_I}$.
\end{lem}
\begin{proof}
The map $\sigma_I$ is an isomorphism. Therefore its restriction to
$\cD (Y)^G$ is one-to-one. Consider the $\bbZ^n$-grading of
$\cD(Y)^G$ given by ~\eqref{AG}. Since
\begin{equation*}
\sigma_I (\Pi_i)=\left\{\begin{array}{l}-\Pi_i -1\mbox{ if }i\in I\\
\Pi_i\mbox{ if }i\notin I
\end{array}\right. ,
\end{equation*}
we can easily see that $\sigma_I (\cA_0)=\cA_0$.
We can also check that
$\sigma_I(u_{\alpha})=\pm u_{\alpha^I}$
with $\alpha^I=(\alpha_1^I,\ldots ,\alpha_s^I)$ where
\begin{equation}\label{alphaI}
\alpha^I_i=\left\{\begin{array}{l}-\alpha_i\mbox{ if }i\in I\\
\alpha_i\mbox{ if }i\notin I\end{array}\right. .
\end{equation}
From this we conclude that $\sigma_I (\cD(Y)_{\alpha
})=\cD(Y)_{\alpha^I}$. Therefore $\sigma_I
(\cD(Y)^G)=\cD(Y)^{G_I}$.
\end{proof}

The automorphism $\sigma_I$ can be thought of as a partial Fourier
transform \cite{C}, Chapter 5. We consider the $\cD (Y)^G$-module
$\cO (Y)^I$ which equals $\cO (Y)$ as a vector space with a new
action given by
\begin{equation}
a\centerdot u=\sigma_I (a) u
\end{equation}
for $a\in\cD (Y)^G$, $u\in\cO (Y)$. Since $\sigma (\cD (Y)^G)=\cD
(Y)^{G_I}$ it follows from Lemma ~\ref{isomorphism} that
\begin{equation}
\cO (Y)^I=\oplus_{\chi\in \Lambda_I}\cO (Y)_{\chi}^I
\end{equation}
where $\Lambda_I=\{L_I\alpha |\alpha\in\bbN^r\times\bbZ^s\}$ and
\begin{equation}
\cO (Y)_{\chi}^I=span\{Q^{\alpha} |L_I \alpha =\chi\}.
\end{equation}
Note in particular that
\begin{equation}
\cO (Y)_0^I=\cO (Y)^{G_I}.
\end{equation}

\section{Main results}\label{ifgw}

The following Lemma can be easily proved.
\begin{lem}\label{hom}
Suppose $R$ is a $\bbZ^s$-graded ring. There is a homomorphism
\begin{displaymath}
\Phi : R\longrightarrow R[x_1^{\pm 1},\ldots ,x_s^{\pm 1}]
\end{displaymath}
defined by $\Phi (r)=rx^{\alpha }$ for $r\in R(\alpha )$,
$\alpha\in\bbZ^s$.
\end{lem}

Assume that $\eta_{r+1},\ldots ,\eta_n$ are linearly independent.
Recall that $G$ acts on $Y$ by the matrix $L$ as in ~\eqref{newL}.

Set
\begin{align}
T_1&=\{\beta\in\bbZ^{r} | L_1\beta=0\}\\
T_2&=\{\beta\in\bbZ^n | L\beta=0\}
\end{align}
and
\begin{align}
T_1'=\{\beta\in T_1| L_2\beta\in d_1\bbZ\times\ldots\times
d_s\bbZ\}.
\end{align}
For $\beta\in T_1'$ let $\kappa (\beta )\in\bbZ^s$ be the column
vector with $i^{th}$ entry $-\gamma_i/d_i$ where
$\gamma=L_2\beta\in\bbZ^s$, and set $\epsilon
(\beta)=\left(\begin{array}{c}\beta\\\kappa
(\beta)\end{array}\right)$.

Part of the Luna slice theorem states that there is a closed
$H$-stable subvariety $S$ containing $w$ and \' etale maps $V/\!/
H\longleftarrow S/\!/ H\longrightarrow Y/\!/ G$, \cite{Lu},
\cite{Sl}. Taking $S=V+w$ we can strengthen this statement in our
case as follows.

\begin{thm}
There is an isomorphism between $V/\!/H$ and $Y/\!/G$.
\end{thm}
\begin{proof} Since $\cO (V)^H=span\{Q^{\beta}
|L_1\beta=0,\;\beta\in\bbN^{r}\cap T_1'\}$ and $\cO
(Y)^G=span\{Q^{\beta}|\beta\in (\bbN^{r}\times\bbZ^{s})\cap T_2\}$
there is an isomorphism between $\cO (V)^H\cong\cO(Y)^G$ sending
$Q^{\beta}$ to $Q^{\epsilon (\beta)}$.
\end{proof}

\subsection{Proof of Proposition C.}
By \cite{MV}, Theorem 8.2.1 and the remarks in Section ~\ref{rem}
we have
\begin{equation}\label{gkin}
\mbox{GKdim}(\cB_{\mu}(V))\leq\mbox{GKdim}(\cB_{\lambda}(Y))=2(n-m).
\end{equation}
 Using the notation of \S ~\ref{notation} we have
\begin{equation*}
\cR=\cD (V)^{H}=\oplus_{\alpha\in T_1'}\cA_0' u_{\alpha}
\end{equation*}
and
\begin{equation*}
\cD (Y)^G=\oplus_{\alpha\in T_2} \cA_0 u_{\alpha}.
\end{equation*}
We regard $\cR$ as a $\bbZ^{s}$-graded ring with deg$(\cA_0'
u_{\alpha})=\kappa (\alpha)$, and apply Lemma ~\ref{hom} to obtain
a homomorphism
\begin{equation}
\cR\longrightarrow \cR[x_1^{\pm 1},\ldots ,x_{s}^{\pm 1}]
\end{equation}
which is the identity on $\cR$ and maps $u_{\alpha}$ to
$u_{\alpha}x^{\kappa (\alpha)}$. Composition with the map
specializing $x_i^{\pm 1}$ to $Q_{r+i}^{\pm 1}$ for $1\leq i\leq
s$ gives the map $\xi$. Note that $\xi (u_{\alpha})=u_{\epsilon
(\alpha )}$.

Since $(\FRAK{h}-\mu(\FRAK{h}))$ is in the kernel of the map from
$\cD (V)^{H}$ to $\cB_{\lambda}(Y)$, $\xi$ induces a map
$\xi_{\lambda}:\cB_{\mu}(V)\longrightarrow \cB_{\lambda}(Y)$. To
show that $\xi_{\lambda }$ is surjective it remains to show that
the image of $\cA_0$ in $\cB_{\lambda}(Y)$ is contained in the
image of $\xi_{\lambda}$. For $1\leq i\leq s$ we have
$d_i\Pi_{r+i}+\sum_{j=1}^r l_{m-s+i,j}\Pi_j\in\FRAK{g}$. Hence
mod$(\FRAK{g}-\lambda(\FRAK{g}))$, $\Pi_{r+i}$ is in the span of
$\Pi_1,\ldots ,\Pi_r$ and $1$. On the other hand
$\xi(\Pi_i)=\Pi_i$ for $1\leq i\leq r$.

Surjectivity of $\xi_{\lambda}$ implies that we have equality in
~\eqref{gkin}. Since $\cB_{\mu}(V)$ is a domain it follows from
\cite{KL}, Proposition 3.15 that $\xi_{\lambda}$ is injective.

\subsection{Proof of Proposition A}

By Lemma \ref{li} it suffices to show that $\cD (Y)^{G}$ has no
finite dimensional modules if $\eta_{r+1},\ldots ,\eta_n$ are
linearly dependent. This follows from \cite{MV}, Proposition
10.1.1(1), but we can give a direct proof as follows. There exist
integers $c_{r+1},\ldots ,c_n$ not all zero such that
$\sum_{i=r+1}^n c_i\eta_i=0$. We can assume that $c = c_n\neq 0$.
Then $\mathcal{Q} = \prod_{i=r+1}^n Q_i^{c_i}$ and $\mathcal{P} =
P_n^c \prod_{i=r+1}^{n-1} Q_i^{-c_i}$ belong to $\cD (Y)^{G}$

Let $\omega = e^{\frac{2\pi i}{c}}$ and consider the automorphism
of $A_1 = k[P,Q]$ sending $P$ to $\omega P$ and $Q$ to
$\omega^{-1}Q$.  Let ${\bbF}$ be the subgroup of $Aut(A_1)$
generated by this automorphism.  We have $A_1^{\bbF} =
k[P^c,Q^c,PQ]$, and it is well known that this is a simple ring,
\cite{Alev}.  Note that $A_1^{\bbF}$ is $\bbZ$-graded when we set
$deg(Q^c) = 1,deg(P^c) = -1$, and $deg(PQ) = 0$.  Applying Lemma
\ref{hom} we see that there is a ring homomorphism $A_1^{\bbF}
\longrightarrow \cD (Y)^{G}$ sending $Q^c$ to $\mathcal{Q}$, $P^c$
to $\mathcal{P}$ and $PQ$ to  $P_nQ_n$. Since $A_1^{\bbF}$ is
simple it follows that this map is injective and $\cD (Y)^{G}$ has
no finite dimensional modules.

\subsection{Proof of Theorem B}
(1)$\implies$(4) If $V^{H^o}=0$ then by Lemma ~\ref{I} there
exists $I\subseteq\{1,\ldots ,r\}$ such that $\cO (Y)^{G_I}=k$.
Therefore by Lemmas ~\ref{special} and ~\ref{isomorphism} the
faithful $\cD (Y)^G$-module $\cO (Y)^I$ is a direct sum of finite
dimensional simple modules $\cO (Y)_{\chi}^I$. Hence $\cD (Y)^G$
has enough finite dimensional simple modules.

(3)$\implies$(1) If $V^{H^o}\neq 0$ we claim that $\cD (Y)^G$ has
no finite dimensional representations. We can assume that
$\rho_1=0$. Then $A_1\cong k[Q_1, P_1]$ is a subalgebra of $\cD
(V)^{H^o}$. Hence the invariants in $A_1$ under the action of the
finite group $H/H^o$ form a subalgebra of $\cD (V)^H$. Since $\cD
(V)^{H}$ is a subalgebra of $\cD (Y)^G$ by Proposition C the claim
follows.

Since (4) obviously implies (3) this proves the equivalence of
(1), (3), (4). The equivalence of (1) and (2) is a special case of
the equivalence of (1) and (4).

\section{More on the modules $\cO (Y)^I$}\label{motm}
We give some alternative descriptions of the modules $\cO (Y)^I$.
First it is easy to see that as $\cD (Y)$-modules
\begin{equation}
\cO (Y)^I\cong \cD (Y)/(\sum_{i\in I}\cD (Y)Q_i +\sum_{i\notin
I}\cD (Y)P_i).
\end{equation}
Next for $I\subseteq\{1,\ldots , r\}$ set $Q_I=\prod_{i\in I} Q_i$
and let
\begin{equation}
\cO (Y)_I=\cO (Y)[Q_I^{-1}].
\end{equation}
Note that $\cO (Y)_I$ is a $\cD (Y)$-module. Now let $M_I$ be the
sum of all submodules of $\cO (Y)_I$ of the form $\cO (Y)_J$ where
$J$ is a proper subset of $I$. Then set
\begin{equation}
N_I=\cO (Y)_I/M_I.
\end{equation}
\begin{prop}
As a $\cD (Y)$-module (and hence also as a $\cD (Y)^G$-module)
$\cO (Y)^I\cong N_I$.
\end{prop}
\begin{proof}
Let $n_I$ be the image of $Q^{-1}_I$ in $N_I$. It is easy to
verify that $Q_i\cdot n_I=0$ if $i\in I$ and $P_i\cdot n_I=0$ if
$i\notin I$. Since $N_I$ is generated by $n_I$ it follows that
there is a surjective map from $\cO (Y)^I$ onto $N_I$. This map is
an isomorphism since $\cO (Y)^I$ is a simple $\cD (Y)$-module.
\end{proof}

\section{FCR-algebras}
\label{FCR}

We give examples of Noetherian FCR-algebras of every given integer
GK dimension $\geq 3$. It is conjectured in \cite{KiS} that FCR
algebras with GK dimension $2$ do not exist. It was proved in
\cite{MV}, Theorem 8.2.1(4) that GK-dim$\cD (Y)^G=2n-m$.

\begin{enumerate}

\item {\bf Odd GK dimension $\geq 3$.} Assume $dimG=1$. We can
take $n\geq 2$. Let the action of $G$ on the $n^{th}$ Weyl algebra
$A_n$ be given by the matrix:
\begin{equation}
\left[
\begin{array}{cccc}
b_{1}&\ldots &b_{n}
\end{array}
\right]
\end{equation}
with $b_{i}\neq 0$ for all $i\in\{1,\ldots ,n\}$. By \cite{MV},
Proposition 10.2.1(3), these are algebras with the reductive
property, and by Theorem B they have enough simple finite
dimensional modules. Then $\mbox{GKdim}A_n^G=2n-1$. In this way we
have many examples of FCR-algebras with odd GK-dim$\geq 3$.

If we take all the $b_i$ equal to $1$ so that $k^\times$ acts by
scalar multiplication on $k^n$ we obtain an unpublished example of
the first author and M. Van den Bergh which is  also mentioned in
\cite{KiS}, \S 3.

\item {\bf Even GK dimension $\geq 6$.} Assume $dim G=2$ and
$n\geq 4$ . Then $\mbox{GK-dim}A_n^G=2n-2$. Let the action of $G$
on $A_{n}$ be given by the matrix:

\begin{equation}
\left[
\begin{array}{ccccccc}
1&0&1&0&\ldots &1&0\\
0&1&0&1&\ldots &0&1
\end{array}
\right]\mbox{ or } \left[
\begin{array}{ccccccc}
1&0&1&0&\ldots &0&1\\
0&1&0&1&\ldots &1&0
\end{array}
\right].
\end{equation}

In this way we get examples of FCR-algebras with even GK-dim$\geq
6$. By  \cite{MV}, Corollary 10.1.6 , $A_n^G$ has the reductive
property, and by Theorem B it has enough simple finite dimensional
modules.

\item {\bf GK-dimension 4.} Take dim$G=2$, $n=3$, $s=1$ and $r=2$.
Consider the algebra
\begin{equation*}
A=k[Q_1,Q_2,Q_3^{\pm 1},P_1, P_2, P_3].
\end{equation*}
Let the action of $G$ on $A$ be given by the matrix
\begin{equation*}
\left[
\begin{array}{ccc}
1&1&0\\
0&0&1
\end{array}
\right].
\end{equation*}
Then GK-dim$A^G=4$. By  \cite{MV}, Corollary 10.1.6 and Theorem B,
$A^G$ is an FCR-algebra.
\end{enumerate}



\begin{thebibliography}{999}

\bibitem[1]{Alev} J. Alev,  Action de groupes sur $A\sb 1( C)$, Ring theory
(Antwerp, 1985), {\em Lecture Notes in Math.} {\bf 1197},
(Springer, Berlin, 1986), 1-9.

\bibitem[2]{C} S.C. Coutinho, {\em A primer of algebraic D-modules}
(Cambridge University Press, 1995).

\bibitem[3]{F} W. Fulton, {\em Introduction to toric varieties}
(Princeton University Press, 1993).

\bibitem[4]{KiS} E.E. Kirkman and L.W. Small, Examples of
FCR-algebras, {\em Comm. Algebra} {\bf 30} (2002), no 7,
3311-3326.

\bibitem[5]{K} H. Kraft, {\em Geometrische Methoden in der Invariantentheorie},
 (Vieweg-Verlag, Braunschweig, 1984).

\bibitem[6]{KS} H. Kraft and L.W. Small, Invariant algebras and
Completely reducible representations, {\em Math. Research Lett.} {\bf 1} (1994),
297-307.

\bibitem[7]{KSW} H. Kraft, L.W.Small and N.R. Wallach, Properties and examples
of FCR-algebras. {\em Manuscripta Math} {\bf 104} (2001), no. 4,
443-450.

\bibitem[8]{KL} G.R. Krause and T.H. Lenagan, Growth of algebras
and Gelfand-Kirillov dimension, {\em Research notes in
mathematics} vol. 116, (Pitman, Boston, 1985).

\bibitem[9]{LS}T. Levasseur and J.T. Stafford, Rings of differential operators
on classical rings of invariants. {\em Mem. Amer. Math. Soc.} {\bf
81} (1989), no. 412.

\bibitem[10]{Lu} D. Luna, Slices \' etales. Sur les groupes
alg\' ebriques, {\em Bull. Soc. Math. France} {\bf 33} (Soc. Math.
France, Paris, 1973), 81-105.

\bibitem[11]{M2} I.M. Musson, Rings of differential operators on
invariant rings of tori, {\em Trans. Amer. Math. Soc.} {\bf 303}
(1987), 805-827.

\bibitem[12]{M} I.M. Musson, Differential operators on toric
varieties,
{\em J. Pure and Applied Algebra} {\bf 95} (1994), 303-315.

\bibitem[13]{MV} I.M. Musson and M. Van den Bergh, Invariants under tori
of rings of differential operators and related topics, {\em Mem.
Amer. Math. Soc.} {\bf 136} (1998), (650).

\bibitem[14]{SG} G.W. Schwarz, Lifting differential operators from
orbit spaces. {\em Ann. Sci. \' Ecole Norm. Sup.} {\bf 4} 28
(1995), no. 3, 253-305.

\bibitem[15]{Sch} G.W. Schwarz, Finite dimensional representations of invariant
differential operators, {\em J. of Alg} {\bf 258} (2002), 160-204.

\bibitem[16]{Sl} P. Slodowy, Der Scheibensatz f$\ddot{u}$r algebraische
Transformationsgruppen (pp. 89--113); Algebraische
Transformationsgruppen und Invariantentheorie. Edited by H. Kraft,
P. Slodowy and T. A. Springer. DMV Seminar, {\bf 13}. Birkhäuser
Verlag, Basel, 1989.

\bibitem[17]{VB} M. Van den Bergh, Some rings of differential operators for
${\rm Sl}\sb2$-invariants are simple. Contact Franco-Belge en
Alg\` {e}bre (Diepenbeek, 1993). J. Pure Appl. Algebra {\bf 107}
(1996), no. 2-3, 309-335.

\end{thebibliography}
\end{document}